\documentclass{amsart}%
\usepackage{amsmath}
\usepackage{graphicx}
\usepackage{amsfonts}
\usepackage{amssymb}%
\newtheorem{theorem}{Theorem}[section]
\theoremstyle{plain}

\newtheorem{proposition}[theorem]{Proposition}

\newtheorem{definition}[theorem]{Definition}
\newtheorem{corollary}[theorem]{Corollary}

\numberwithin{equation}{section}
\begin{document}
\title[Higher order Nielsen numbers.]{Higher order Nielsen numbers.}
\author{Peter Saveliev}
\address{Marshall University, Huntington, WV 25755}
\email{saveliev@member.ams.org}
\urladdr{http://www.saveliev.net}
\subjclass{Primary 55M20, 55H25; Secondary 55N22.}
\keywords{Nielsen theory, coincidence index, preimage problem, bordism}

\begin{abstract}
Suppose $X,Y$ are manifolds, $f,g:X\rightarrow Y$ are maps. The well-known
Coincidence Problem studies the coincidence set $C=\{x:f(x)=g(x)\}.$ The
number $m=\dim X-\dim Y$ is called the codimension of the problem. More
general is the Preimage Problem. For a map $f:X\rightarrow Z$ and a
submanifold $Y$ of $Z,$ it studies the preimage set $C=\{x:f(x)\in Y\},$ and
the codimension is $m=\dim X+\dim Y-\dim Z.$ In case of codimension $0$, the
classical Nielsen number $N(f,Y)$ is a lower estimate of the number of points
in $C$ changing under homotopies of $f,$ and for an arbitrary codimension, of
the number of components of $C$. We extend this theory to take into account
other topological characteristics of $C.$ \textit{The goal is to find a
\textquotedblleft lower estimate\textquotedblright\ of the bordism group
}$\Omega_{p}(C)$\textit{ of }$C.$ The answer is the Nielsen group $S_{p}(f,Y)$
defined as follows. In the classical definition the Nielsen equivalence of
points of $C$ based on paths is replaced with an equivalence of singular
submanifolds of $C$ based on bordisms. We let $S_{p}^{\prime}(f,Y)=\Omega
_{p}(C)/\sim_{N},$ then the Nielsen group of order $p$ is the part of
$S_{p}^{\prime}(f,Y)$ preserved under homotopies of $f$. The Nielsen number
$N_{p}(F,Y)$ of order $p$ is the rank of this group (then $N(f,Y)=N_{0}%
(f,Y))$. These numbers are new obstructions to removability of coincidences
and preimages. Some examples and computations are provided.

\end{abstract}
\maketitle

\section{Introduction}

Suppose\ $X,Y$ are smooth orientable compact manifolds, $\dim X=n+m,$ $\dim
Y=n$, $m\geq0$ the \textit{codimension}, $f,g:X\rightarrow Y$ are maps, the
\textit{coincidence set}
\[
C=Coin(f,g)=\{x\in X:f(x)=g(x)\}
\]
is a compact subset of $X\backslash\partial X.$

Consider the Coincidence Problem: \textquotedblleft What can be said about the
coincidence set $C$ of $(f,g)$?\textquotedblright\ One of the main tools is
the Lefschetz number $L(f,g)$ defined as the alternating sum of traces of a
certain endomorphism on the homology group of $Y$. The famous Lefschetz
coincidence theorem provides a sufficient condition for the existence of
coincidences (codimension $m=0$): $L(f,g)\neq0\Longrightarrow C=Coin(f,g)\neq
\emptyset$, see \cite[VI.14]{Bredon}, \cite[Chapter 7]{Vick}$.$

Now, what else can be said about the coincidence set? As $C$ changes under
homotopies of $f$ and $g,$ a reasonable approach is to \ try to minimize the
\textquotedblleft size\textquotedblright of $C.$ In case of zero codimension
$C$ is discrete and we simply minimize the number of points in $C$. The result
is the Nielsen number. It is defined as follows. Two points $p$,$q\in C$
belong to the same Nielsen class if (1) there is a path $s$ in $X$ between $p$
and $q$; (2) $fs$ and $gs$ are homotopic relative to the end-points$.$ A
Nielsen class is called essential if it cannot be removed by a homotopy of
$f,g$ (alternatively, a Nielsen class is algebraically essential if its
coincidence index is nonzero \cite{BB}). Then the Nielsen number $N(f,g)$ is
the number of essential Nielsen classes. It is a lower estimate of the number
of points in $C$. In case of positive codimension $N(f,g)$ still makes sense
as a lower estimate of the number of components of $C$ \cite{Wong}. However
only for $m=0$ the Nielsen number is known to be a \textit{sharp} estimate,
i.e., there are maps $f^{\prime},g^{\prime}$ compactly homotopic of $f,g$ such
that $C^{\prime}=Coin(f^{\prime},g^{\prime})$ consists of exactly $N(f,g)$
path components (Wecken property). This minimization is achieved by removing
inessential classes through homotopies of $f,g.$

The Nielsen theory for codimension $m=0$ is well developed, for the Fixed
Point and the Root Problems \cite{Brown}, \cite{Jiang}, \cite{Kiang}, and for
the Coincidence Problem \cite{BS}. However, for $m>0,$ the vanishing of the
coincidence index does not guarantee that the Nielsen class can be removed.
Some progress has been made for codimension $m=$ $1.$ In this case the
secondary obstruction to the removability of a coincidence set was considered
by Fuller \cite{Fuller} for $Y$ simply connected. Hatcher and Quinn \cite{HQ}
showed that the obstruction to a higher dimensional Whitney Lemma lies in a
certain framed bordism group. Based on this result, necessary and sufficient
conditions of the removability of a Nielsen class were studied by Dimovski and
Geoghegan \cite{DG}, Dimovski \cite{Dim} for parametrized fixed point theory,
i.e., when $f:Y\times I\rightarrow Y$ is the projection. The results of
\cite{DG} were generalized by Jezierski \cite{Jez} for the coincidence problem
$f,g:X\rightarrow Y,$ where $X,Y$ are open subsets of Euclidean spaces or $Y$
is parallelizable. Geoghegan and Nicas \cite{GN} developed a parametrized
Nielsen theory based on Hochschild homology. For some $m>1,$ sufficient
conditions of the local removability are provided in \cite{Sav2}. Necessary
conditions of the global removability for arbitrary codimension are considered
by Gon\c{c}alves, Jezierski, and Wong \cite[Section 5]{GJW} with $N$ a torus
and $M$ a nilmanifold.

In these papers higher order Nielsen numbers are not explicitly defined
(except for \cite{Dim}, see the comment in the end of the paper). However they
all contribute to the problem of finding the lower estimate of the number of
components of $C$. We extend these results to take into account other
topological characteristics of $C.$ In the spirit of the classical Nielsen
theory, \textit{our goal is to find \textquotedblleft lower
estimates\textquotedblright\ of the bordism groups }$\Omega_{\ast}(C).$

The crucial motivation for our approach is the removability results for
codimension 1 due to Dimovski and Geoghegan \cite{DG} and Jezierski
\cite{Jez}. Consider Theorem 5.3 in \cite{Jez}. Assume that codimension $m=$
$1,$ $n\geq4,$ $X,Y$ are open subsets of Euclidean spaces. Suppose $A$ is a
Nielsen class. Then if $f,g$ are transversal, $A$ is the union of disjoint
circles. Define the Pontriagin-Thom map (PT) as the composition
\[
\mathbf{S}^{n+1}\simeq\mathbf{R}^{n+1}\cup\{\ast\}\rightarrow\mathbf{R}%
^{n+1}/(\mathbf{R}^{n+1}\backslash\nu)\simeq\overline{\nu}/\partial
\nu^{\underrightarrow{~\ \ \ f-g~~}}\mathbf{R}^{n}/(\mathbf{R}^{n}\backslash
D)\simeq\mathbf{S}^{n},
\]
where $\nu$ is a normal bundle of $A,$ $D\subset\mathbf{R}^{n}$ is a ball
centered at $0$ satisfying $(f-g)(\partial\nu)\subset\mathbf{R}^{n}\backslash
D.$ It is an element of $\pi_{n}(\mathbf{S}^{n-1})=\mathbf{Z}_{2}.$ Then $A$
can be removed if and only if the following conditions are satisfied:

$\text{(W1) }A=\partial S,\text{ }$\textit{where} $S$ \textit{is an orientable
connected surface}$,\text{ }f|_{S}\sim g|_{S}\text{ rel }A$\textit{ (the
surface condition)}$;$

(W2)\textit{ the PT map is trivial (the }$\mathbf{Z}_{2}$\textit{-condition).}%
$\ \ \ \ $

Earlier Dimovski and Geoghegan \cite{DG} considered a similar pair of
conditions (not independent though) in their Theorem 1.1 and compared them to
the codimension $0$ case. They write: \textquotedblleft...the role of `being
in the same fixed point class' is played by the surface condition (i), while
that of the fixed point index is played by the natural orientation. The
$\mathbf{Z}_{2}$-obstruction is a new feature...\textquotedblright\ One can
use the first observation to define the Nielsen equivalence on the set of of
$1$-submanifolds of $C$ (here $A$ is Nielsen equivalent to the empty set).
However, we will see that the PT map has to serve as the index of the Nielsen
class. The index will be defined in the traditional way but with respect to an
arbitrary homology theory $h_{\ast}.$ Indeed in the above situation it is an
element of the stable homotopy group $\pi_{n+1}^{S}(\mathbf{S}^{n}%
)=\mathbf{Z}_{2}$.

More generally, we define the Nielsen equivalence on the set $M_{m}(C)$ of all
closed singular $m$-manifolds in $C=Coin(f,g).$ Two singular $m$-manifolds
$p:P\rightarrow C$ and $q:Q\rightarrow C$ belong to the same \textit{Nielsen
class}, $p\sim_{N}q,$ if

1. $ip$ and $iq$ are bordant, where $i:C\rightarrow N$ is the inclusion, i.e.
there is a map $F:W\rightarrow N$ extending $ip\sqcup iq$ such that $W$ is a
bordism between $P$ and $Q;$

2. $fF$ and $gF$ are homotopic relative to $fp,fq.$

\noindent Then $S_{m}^{\prime}(f,g)=M_{m}(C)/\sim_{N}$ is the group of Nielsen
classes. Let $S_{m}^{a}(f,g)$ be the group of algebraically essential Nielsen
classes, i.e., the ones with  non-trivial index. Then the (algebraic) Nielsen
number of order $m$ is the rank of $S_{m}^{a}(f,g)$ (these numbers are new
obstructions to removability of coincidences). In light of this definition
Jezierski's theorem can be thought of as a Wecken type theorem for $m=1$.

An area of possible applications of the coincidence theory for positive
codimension is discrete dynamical systems (for a related theory of flows see
\cite{GN}). A \textit{dynamical system} on a manifold $M$ is determined by a
map $f:M\rightarrow M.$ Then the next position $f(x)$ depends only on the
current one, $x\in M$. Suppose now that we have a fiber bundle $F\rightarrow
N^{\underrightarrow{~\ \ \ g\ \ \ ~~}}M$ and a map $f:N\rightarrow M.$ This is
a \textit{parametrized\ dynamical system}, where the next position $f(x,s)$
depends not only on the current one, $x\in M,$ but also the current
\textquotedblleft state\textquotedblright, $s\in F.$ The coincidence theory
describes the set of all positions and states such that the positions remain
fixed, $f(x,s)=x$. Alternatively, $x$ is the \textquotedblleft
state\textquotedblright\ and $s$ is the \textquotedblleft
input\textquotedblright\ in a control system. A \textit{control system}
\cite[p. 16]{Nij} is defined as a commutative diagram
\[%
\begin{tabular}
[c]{ccc}%
$N$ & $^{\underrightarrow{\quad h\quad}}$ & $TM,$\\
$\downarrow^{\pi}$ & $\swarrow^{\pi_{M}}$ & \\
$M$ &  &
\end{tabular}
\ \
\]
where $N$ is a fiber bundle over $M.$ The translation along trajectories of
this system creates a parametrized dynamical system.

Instead of the Coincidence Problem, throughout the rest of the paper we apply
the approach outlined above to the Nielsen Theory for the so-called
\textit{Preimage Problem }considered by Dobrenko and Kucharski \cite{DK}.
Suppose\ $X,Y,Z$ are connected CW-complexes, $Y\subset Z,$ $f:X\rightarrow Z$
is a map. The problem studies the set $C=f^{-1}(Y)$ and can be easily
specialized to the Fixed Point Problem if we put $Z=X\times X,$ $Y=d(X),$
$f=(Id,g),$ to the Root Problem if $Y$ is a point$,$ and to the Coincidence
Problem if $Z=Y\times Y,$ $Y$ is the diagonal of $Z,$ $f=(F,G)$ (see
\cite{McCord})$.$

Suppose $X,Y,Z$ are smooth manifolds and $f$ is transversal to $Y.$ Then under
the restriction $\dim X+\dim Y=\dim Z,$ the preimage $C=f^{-1}(Y)$ of $Y$
under $f$ is discrete$.$ The Nielsen number $N(f,Y)$ is the sharp lower
estimate of the least number of points in $g^{-1}(Y)$ for all maps $g$
homotopic to $f$ \cite[Theorem 3.4]{DK}$,$ i.e., $N(f,Y)\leq\#g^{-1}(Y)$ for
all $g\sim f.$ If we omit the above restriction, $C$ is an $r$-manifold
\cite[Theorem II.15.2, p. 114]{Bredon}, where
\[
r=\dim X+\dim Y-\dim Z.
\]

\textit{The Setup. }$X,Y,Z$ are connected CW-complexes, $Y\subset Z,$
\[
\dim X=n+m,\dim Y=n,\dim Z=n+k,
\]
$f:X\rightarrow Z$ is a map, the \textit{preimage set} $C=f^{-1}(Y),$ the
\textit{codimension} of the problem
\[
r=n+m-k,
\]
$j:C\rightarrow X$ the inclusion.

The paper is organized as follows. Just as for the coincidence problem we
define the Nielsen equivalence of singular $q$-manifolds in $C$ and the group
of Nielsen classes $S_{q}^{\prime}(f)=M_{q}(C)/\sim_{N}=\Omega_{q}(C)/\sim
_{N},$ where $\Omega_{\ast}$ is the orientable bordism group (Section
\ref{NClass})$.$ Next we identify the part of $S_{q}^{\prime}(f)$ preserved
under homotopies of $f.$ The result is the Nielsen group $S_{q}(f)$, the group
of topologically essential classes (Section \ref{TopEss})$.$ As we have
described above, \textit{the Nielsen group is a subgroup of a quotient group
of }$\Omega_{q}(C)$\textit{ and, in this sense, its \textquotedblleft lower
estimate\textquotedblright.}

\begin{proposition}
$S_{\ast}(f)$ is homotopy invariant.
\end{proposition}

The Nielsen number of order $p$\textit{, }$p=0,1,2,...,$ is defined as
$N_{p}(f,Y)=\operatorname*{rank}\,S_{p}(f,Y).$ Clearly the classical Nielsen
number is equal to $N_{0}(f).$

\begin{proposition}
$N_{p}(f)\leq\operatorname*{rank}\,\Omega_{p}(g^{-1}(Y))$ if $f\sim g.$
\end{proposition}

In Section \ref{Natur} we discuss the naturality of the Nielsen group. In
particular we obtain the following

\begin{proposition}
Given $Z,Y\subset Z.$ Then $S_{\ast}$ is a functor from the category of
preimage problems as pairs $(X,f),$ $f:X\rightarrow Z,$ with morphisms as maps
$k:X\rightarrow U$ satisfying $gk=f,$ to the category of graded abelian groups.
\end{proposition}

For the manifold case, there is an alternative approach to essentiality. In
Section \ref{BordInd} the \textquotedblleft preimage index\textquotedblright%
\ is defined simply as $I_{f}=f_{\ast}:$ $\Omega_{\ast}(C)\rightarrow
\Omega_{\ast}(Y).$ It is a homomorphism on $S_{\ast}^{\prime}(f)$ and the
group of algebraically essential Nielsen classes is defined as $S_{\ast}%
^{a}(f,Y)=S_{\ast}^{\prime}(f,Y)/\ker I_{f}.$ We show that every algebraically
essential class is topologically essential. In Section \ref{Index} we consider
the traditional index $Ind_{f}(P)$ of an isolated subset $P$ of $C$ in terms
of a generalized homology $h_{\ast}.$ It is defined in the usual way as the
composition
\[
h_{\ast}(X,X\backslash U)^{\underleftarrow{\quad\simeq\quad\ }}h_{\ast
}(V,V\backslash U)^{\underrightarrow{\quad f_{\ast}\quad}}h_{\ast
}(Z,Z\backslash Y),
\]
where $V\subset\overline{V}\subset U$ are neighborhoods of $P.$ Then we show
how it is related to $I_{f}.$ In Section \ref{Examples} we consider some
examples of computations of these groups, especially in the setting of the
Pontriagin-Thom construction. 

In Sections \ref{Wecken} and \ref{Jez} based on Jezierski's theorem we prove
the following Wecken type theorem for codimension $1.$

\begin{proposition}
Under conditions of Jezierski's theorem, $f,g$ is homotopic to $f^{\prime
},g^{\prime}$ such that
\[
S_{p}(f,g)\simeq\Omega_{p}(Coin(f^{\prime},g^{\prime})),p=0,1.
\]

\end{proposition}

This approach to the study of the homology of the preimage set has a potential
application in robust control \cite{Jon}. The so-called Nyquist map is a map
from the \textquotedblleft uncertainty space\textquotedblright\ (i.e., the
state input manifold) of the control system to the Riemann sphere. Topological
robust control theory tracks the change of the topology of the preimage of the
origin under perturbations of this map. In general, robust stability criteria
consider \textquotedblleft mapping of the uncertainty into a `performance
evaluation' space... and checking whether the image is in the correct
subset\textquotedblright\ \cite[p. 20]{Jon}.

To motivate our definitions, in the beginning of each section be will review
the relevant part of Nielsen theory for the Preimage Problem following
Dobrenko and Kucharski \cite{DK}, and McCord \cite{McCord}.

All manifolds are assumed to be orientable and compact.

\section{Nielsen Classes.\label{NClass}}

In Nielsen theory, two points $x_{0},x_{1}\in C=f^{-1}(Y)$ belong to the same
Nielsen class, $x_{0}\sim x_{1},$ if

(1) there is a path $\alpha:I\rightarrow X$ such that $\alpha(i)=x_{i};$

(2) there is a path $\beta:I\rightarrow Y$ such that $\beta(i)=f(x_{i});$

(3) $f\alpha$ and $\beta$ are homotopic relative to $\{0,1\}.$

This is an equivalence relation partitioning $C$ into a finite number of
Nielsen classes. However, since we want Nielsen classes to form a group, we
should think of $x_{0},x_{1}$ as singular $0$-manifolds in $C$ (a singular
$p$-manifold in $M$ is a map $s:N\rightarrow M,$ where $N$ is a $p$-manifold).
Then conditions (1) and (2) express the fact that $x_{0},x_{1}$ are bordant in
$X,$ and $f(x_{0}),f(x_{1})$ are bordant in $Y.$

Recall \cite{CF}, \cite{Stong} that two closed $p$-manifolds $N_{0},N_{1}$ are
called \textit{bordant} if there is a \textit{bordism} between them, i.e., a
$(p+1)$-manifold $W$ such that $\partial W=N_{0}\sqcup-N_{1}.$ Two closed
singular manifolds $s_{i}:N_{i}\rightarrow M,$ $i=0,1,$ are \textit{bordant},
$s_{0}\sim_{b}s_{1},$ if there is a map $h:W\rightarrow M$ extending
$s_{0}\sqcup s_{1},$ where $W$ is a bordism between $N_{0}$ and $N_{1}.$

Let $M_{p}(A,B)$ denote the set of all singular $p$-manifolds $s:(N,\partial
N)\rightarrow(A,B)$.

\begin{definition}
\label{MainDef}Two singular $p$-manifolds $s_{0},s_{1}\in M_{p}(C)$ in $C,$
i.e., maps $s_{i}:S_{i}\rightarrow C,$ $i=0,1,$ are \textit{Nielsen
equivalent}, $s_{0}\sim_{N}s_{1},$ if

(1) $js_{0},js_{1}$ are bordant in $X$ via a map $H:W\rightarrow X$ extending
$s_{0}\sqcup s_{1}$ such that $W$ is a bordism between $S_{0}$ and $S_{1};$

(2) $fs_{0},fs_{1}$ are bordant in $Y$ via a map $G:W\rightarrow Y$ extending
$fs_{0}\sqcup fs_{1};$

(3) $fH$ and $G$ are homotopic relative to $S_{0}\sqcup S_{1}.$
\end{definition}

We denote the Nielsen class of $s\in M_{p}(C)$ by $[s]_{N},$ or simply $[s].$

\begin{proposition}
$\sim_{N}$ is an equivalence relation on $M_{p}(C).$
\end{proposition}

\begin{definition}
\textit{The group of Nielsen classes of order }$p$, $S_{p}^{\prime}(f,Y),$ or
simply $S_{p}^{\prime}(f),$ is defined as
\[
S_{p}^{\prime}(f,Y)=M_{p}(C)/\sim_{N}.
\]
\textit{The group of Nielsen classes for the coincidence problem will be
denoted by} $S_{p}^{\prime}(f,g).$
\end{definition}

In contrast to the classical Nielsen theory, the elements of Nielsen classes
are not points but sets of points. Even in the case of $p=0,$ one has more to
deal with. For example, suppose $C=\{x,y\}$ and $x\sim_{N}y.$ The the elements
of $S_{0}^{\prime}(f)$ are $[\{x,y\}]=[\{x\}],$
$[\{-x,-y\}]=[\{-x\}]=-[\{x\}],$ $[\{x\}\cup\{-y\}]=[\varnothing],$
$[\{x\}\cup\{y\}]=[\{x\}\cup\{x\}]=[2\{x\}]=2[\{x\}],$ etc.

Another example. Suppose $X=Z=\mathbf{S}^{2},$ $Y$ is the equator of $Z,$ $f$
a map of degree $2$ such that $C=f^{-1}(Y)$ is the union of two circles
$C_{1}$ and $C_{2}$ around the poles. Then $S_{1}^{\prime}(f)=\mathbf{Z}$
generated by $C_{1}\sqcup C_{2}.$ A similar construction applies to
$X=Z=\mathbf{S}^{n},$ $Y=\mathbf{S}^{n-1},$ $n\geq2,$ then $S_{n-1}^{\prime
}(f)=\mathbf{Z}$ is generated by the union of two copies of $\mathbf{S}^{n-1}$.

Let $M_{p}^{h}(A,B)$ denote the semigroup of all homotopy classes, relative to
boundary, of maps $s\in M_{p}(A,B).$ Consider the commutative diagram:
\begin{equation}%
\begin{array}
[c]{ccccc}%
M_{p+1}^{h}(X,C) & ^{\underrightarrow{\qquad\delta\qquad}} & M_{p}^{h}(C) &
^{\underrightarrow{\qquad j_{\ast}\qquad}} & M_{p}^{h}(X)\\
\downarrow^{f_{\ast}} &  & \downarrow^{f_{\ast}} &  & \downarrow^{f_{\ast}}\\
M_{p+1}^{h}(Z,Y) & ^{\underrightarrow{\qquad\delta\qquad}} & M_{p}^{h}(Y) &
^{\underrightarrow{\qquad k_{\ast}\qquad}} & M_{p}^{h}(Z)\\
\uparrow^{I_{\ast}} & \overset{\delta}{\nearrow} &  &  & \\
M_{p+1}^{h}(Y,Y) &  &  &  &
\end{array}
\label{diagram}%
\end{equation}
where $\delta$ is the boundary map, $I$ is the inclusion. Then we have an
alternative way to define the group of Nielsen classes:
\[
S_{p}^{\prime}(f,Y)=M_{p}^{h}(C)/\delta(f_{\ast}^{-1}(\operatorname{Im}%
I_{\ast})).
\]

Let $\Omega_{p}(A,B)$ denote the group of bordism classes in $M_{p}(A,B)$ with
$\sqcup$ as addition. Then $\Omega_{\ast}$ is a generalized homology
\cite{CF}, \cite{Stong}.

\begin{proposition}
If $s_{0}\sim_{N}s_{1}\sim_{b}s_{2}$ then $s_{0}\sim_{N}s_{2}.$ Therefore
$\sim_{N}$ is an equivalence relation on $\Omega_{\ast}(C).$
\end{proposition}

\begin{proposition}
If $s_{0}\sim_{N}s_{1},$ $t_{0}\sim_{N}t_{1}$ then $s_{0}\sqcup t_{0}\sim
_{N}s_{1}\sqcup t_{1}.$ Therefore $\sim_{N}$ is preserved under the operation
of $\Omega_{\ast}(C).$ Thus $S_{\ast}^{\prime}(f,Y)=\Omega_{\ast}(C)/\sim_{N}$
is a group.
\end{proposition}

Next we discuss the naturality of this group.

\begin{definition}
\label{def k*}Suppose we have another preimage problem $f^{\prime}:X^{\prime
}\rightarrow Z^{\prime}\supset Y^{\prime}$ connected to the first by maps
$k:X\rightarrow X^{\prime}$ and $h:Z\rightarrow Z^{\prime}$ such that
$f^{\prime}k=hf$ and $h(Y)\subset Y^{\prime}$ (see diagram in Proposition
\ref{composition})$.$ Then we define the \textit{map induced by }%
$k$\textit{\ and }$h,$%
\[
k_{\ast}^{\prime}:S_{\ast}^{\prime}(f,Y)\rightarrow S_{\ast}^{\prime
}(f^{\prime},Y^{\prime}),
\]
by $k_{\ast}^{\prime}([s]_{N})=[ks]_{N}.$
\end{definition}

\begin{proposition}
\label{k*}$k_{\ast}^{\prime}$ is well defined.
\end{proposition}

\begin{proof}
Let $C^{\prime}=f^{\prime-1}(Y^{\prime}).$ If $x\in C$ then $f(x)=y\in Y.$ Let
$x^{\prime}=k(x)$ and $y^{\prime}=h(y)\in h(Y)\subset Y^{\prime}.$ Then by
assumption $g(x^{\prime})=y^{\prime},$ so $x^{\prime}\in C^{\prime}.$
Therefore the following diagram commutes:
\[%
\begin{array}
[c]{ccc}%
(X,C) & ^{\underrightarrow{\qquad k\qquad}} & (X^{\prime},C^{\prime})\\
\downarrow^{f} &  & \downarrow^{f^{\prime}}\\
(Z,Y) & ^{\underrightarrow{\qquad h\qquad}} & (Z^{\prime},Y^{\prime}).
\end{array}
\]
The second preimage problem has a diagram analogous to (\ref{diagram}).
Together they provide two opposite faces of a 3-dimensional diagram with other
faces supplied by the diagram above. The diagram commutes. Therefore for each
$s\in M_{p}(C),$ $s\sim_{N}\emptyset\Longrightarrow ks\sim_{N}\emptyset.$
\end{proof}

\begin{proposition}
\label{composition}Suppose the following diagram for three preimage problems
commutes:
\[%
\begin{array}
[c]{ccccc}%
Y & ^{\underrightarrow{\quad h\quad}} & Y^{\prime} & ^{\underrightarrow{\quad
l\quad}} & Y^{\prime\prime}\\
\downarrow &  & \downarrow &  & \downarrow\\
Z & ^{\underrightarrow{\quad h\quad}} & Z^{\prime} & ^{\underrightarrow{\quad
l\quad}} & Z^{\prime\prime}\\
\uparrow^{f} &  & \uparrow^{f^{\prime}} &  & \uparrow^{f^{\prime\prime}}\\
X & ^{\underrightarrow{\quad k\quad}} & X^{\prime} & ^{\underrightarrow{\quad
j\quad}} & X^{\prime\prime}%
\end{array}
\]
Then $j_{\ast}^{\prime}k_{\ast}^{\prime}=(jk)_{\ast}^{\prime}:S_{\ast}%
^{\prime}(f,Y)\rightarrow S_{\ast}^{\prime}(f^{\prime\prime},Y^{\prime\prime
}).$
\end{proposition}

\begin{proof}
From the definition, $(jk)_{\ast}^{\prime}([s]_{N})=[jks]_{N}$ and $j_{\ast
}^{\prime}k_{\ast}^{\prime}([s]_{N})=j_{\ast}^{\prime}([ks]_{N})=[jks]_{N}.$
\end{proof}

\begin{proposition}
\label{identity}$(Id_{X})_{\ast}^{\prime}=Id_{S_{\ast}^{\prime}(f,Y)}.$
\end{proposition}

\begin{corollary}
If $\mathcal{P}$ is the category of preimage problems as quadruples
$(X,Z,Y,f),$ $Y\subset Z,$ $f:X\rightarrow Z,$ with morphisms as pairs of maps
$(k,h)$ satisfying Definition \ref{def k*}$,$ then $S_{\ast}^{\prime}$ is a
functor from $\mathcal{P}$ to $\mathbf{Ab}_{\ast},$ the graded abelian groups.
\end{corollary}

\section{Topologically Essential Nielsen Classes.\label{TopEss}}

In the classical theory, a Nielsen class is called essential if it cannot be
removed by a homotopy. More precisely, suppose $F:I\times X\rightarrow Z$ is a
homotopy of $f,$ then the $t$-section $N_{t}=\{x\in X:(t,x)\in N\},$ $0\leq
t\leq1,$ of the Nielsen class $N$ of $F$ is a Nielsen class of $f_{t}%
=F(t,\cdot)$ or is empty \cite[Corollary 1.5]{DK}. Next we say that the
Nielsen classes $N_{0},N_{1}$ of $f_{0},f_{1}$ respectively are in the
$F$-Nielsen relation if there is a Nielsen class $N$ of $F$ such that
$N_{0},N_{1}$ are the $0$- and $1$-sections of $N.$ This establishes an
\textquotedblleft equivalence\textquotedblright\ relation between some Nielsen
classes of $f_{0}$ and some Nielsen classes of $f_{1}.$ Given a Nielsen class
$N_{0}$ of $f_{0},$ if for any homotopy there is a Nielsen class of $f_{1}$
corresponding to $N_{0}$ then $N_{0}$ is called essential. In our theory the
$F$-Nielsen relation takes a simple form of two homomorphisms from $S_{\ast
}^{\prime}(f_{0}),S_{\ast}^{\prime}(f_{1})$ to $S_{\ast}^{\prime}(F).$

Suppose $F:I\times X\rightarrow Z$ is a homotopy, $f_{t}(\cdot)=F(t,\cdot
):X\rightarrow Z,$ and let $i_{t}:X\rightarrow\{t\}\times X\rightarrow I\times
X$ be the inclusions. Since $f_{t}=Fi_{t}$, the homomorphism $i_{t\ast
}^{\prime}:S_{\ast}^{\prime}(f_{t})\rightarrow S_{\ast}^{\prime}(F)$ is well
defined for each $t\in\lbrack0,1]$ (Proposition \ref{k*}). The following
result is crucial.

\begin{theorem}
\label{k-monic}Suppose $F:I\times X\rightarrow Z$ is a homotopy of $f,$
$F|_{\{0\}\times X}=f.$ Suppose $i:X\rightarrow\{0\}\times X\rightarrow
I\times X$ is the inclusion. Then $i_{\ast}^{\prime}:S_{\ast}^{\prime}(f)$
$\rightarrow S_{\ast}^{\prime}(F)$ is injective.
\end{theorem}

\begin{proof}
(cf. \cite[Lemma 1.4]{DK}) Suppose $v\in M_{p}(f^{-1}(Y)),$ $v:M\rightarrow
f^{-1}(Y),$ where $M$ is a $p$-manifold. Then $u=iv=\{0\}\times v\in
M_{p}(F^{-1}(Y)),$ so that $u:M\rightarrow F^{-1}(Y)\supset\{0\}\times
f^{-1}(Y).$ Suppose $[u]_{N}=0$ in $S_{p}^{\prime}(F),$ then there is a $U\in
M_{p+1}(I\times X,F^{-1}(Y)),$ $U:(W,\partial W)\rightarrow(I\times
X,F^{-1}(Y)),$ such that $M=\partial W,$ $U|_{M}=u,$ and $FU:(W,\partial
W)\rightarrow(Z,Y)$ is homotopic relative to $M=\partial W$ to a $G\in
M_{p+1}(Y,Y).$ Then $U=(P,V),$ where $P:W\rightarrow I,$ $P|_{M}=\{0\},$ and
$V:W\rightarrow X,$ $V|_{M}=v.$ Define a homotopy $H:I\times W\rightarrow Z$
by
\[
H(s,x)=F((1-s)P(x),V(x)).
\]
Then $H(0,x)=F(P(x),$ $V(x))=FU(x),$ $H(1,x)=F(0,V(x))=fV(x).$ Suppose $x\in
M.$ Then first, $H(s,x)=F((1-s)\cdot0,v(x))=F(0,v(x))=f(v(x));$ second,
$FU(x)=Fu(x)=Fiv(x)=fv(x);$ third, $fV(x)=fv(x).$ Thus $FU$ and $fV$ are
homotopic relative to $M.$ Therefore $fV$ is homotopic to $G$ relative to $M.$
We have proven that if $[u]_{N}=i_{\ast}^{\prime}[v]_{N}=0$ in $S_{p}^{\prime
}(F)$ then $[v]_{N}=0$ in $S_{p}^{\prime}(f).$ Therefore $\ker i_{\ast
}^{\prime}=\{0\}.$
\end{proof}

Thus the Nielsen classes of a map are included in the Nielsen classes of its homotopy.

This theorem generalizes both the fact that the intersection of a Nielsen
class of $F$ with $\{0\}\times X$ is a Nielsen class of $f_{0}$
\cite[Corollary 1.5]{DK}, for codimension 0, and the fact that (W1) is
homotopy invariant \cite[Lemma 4.2]{Jez}, for codimension 1 (see Section
\ref{Wecken}).

Now the following are monomorphisms
\[
S_{\ast}^{\prime}(f_{0})^{\underrightarrow{\qquad i_{0\ast}^{\prime}\qquad}%
}S_{\ast}^{\prime}(F)^{\underleftarrow{\qquad i_{1\ast}^{\prime}\qquad\ }%
}S_{\ast}^{\prime}(f_{1}).
\]
Let
\[
M_{\ast}^{F}=\operatorname{Im}i_{0\ast}^{\prime}\cap\operatorname{Im}i_{1\ast
}^{\prime}.
\]
Then $M_{\ast}^{F}$ is isomorphic to some subgroups of $S_{\ast}^{\prime}(F),$
$S_{\ast}^{\prime}(f_{0}),$ $S_{\ast}^{\prime}(f_{1})$ (as a subgroup of
$S_{\ast}^{\prime}(f_{0}),$ $M_{\ast}^{F}$ should be understood as the set of
Nielsen classes of $f_{0}$ preserved by $F).$ Now we say that a class
$s_{0}\in S_{\ast}^{\prime}(f_{0})$ of $f_{0}$ is $F$\textit{-related} to a
class $s_{1}\in S_{\ast}^{\prime}(f_{1})$ of $f_{1}$ if there is $s\in
S_{\ast}^{\prime}(F)$ such that $i_{0\ast}^{\prime}(s_{0})=s=i_{1\ast}%
^{\prime}(s_{1}).$ Then $s_{1}=i_{1\ast}^{\prime-1}i_{0\ast}^{\prime}(s_{0})$
if defined, otherwise we can set $s_{1}=0.$ Thus some classes cannot be
reduced to zero by a homotopy and we call them \textit{(topologically)
essential Nielsen classes. }Together (plus zero) they form a group, as follows.

\begin{definition}
\textit{The group of (topologically) essential Nielsen classes} is defined as
\[
S_{\ast}(f,Y)=%
{\displaystyle\bigcap}
\{M_{\ast}^{F}:F\text{ is a homotopy of }f\}\subset S_{\ast}^{\prime}(f,Y).
\]
($S_{p}(f,Y)$ can also be called the \textit{Nielsen group of order }$p$,
while $S_{p}^{\prime}(f,Y)$ the \textit{pre-Nielsen group}.)
\end{definition}

If $f\sim g$ then $S_{\ast}(f)$ $\simeq S_{\ast}(g).$ Therefore,

\begin{theorem}
$S_{\ast}(f)$ is homotopy invariant. Moreover for any $g$ homotopic to $f$
there is a monomorphism $S_{\ast}(f)$ $\rightarrow S_{\ast}^{\prime}(g).$
\end{theorem}

Now $S_{\ast}(f)$ is a subgroup of, which is $S_{p}^{\prime}(f)$ is a quotient
of $\,\Omega_{\ast}(f^{-1}(Y)).$ In this sense, $S_{\ast}(f)$ is a
\textquotedblleft lower estimate\textquotedblright\ of $\,\Omega_{\ast}%
(g^{-1}(Y))$ for any $g$ homotopic to $f$.

\begin{definition}
\textit{The Nielsen number of order }$p$\textit{, }$p=0,1,2,...,$ is defined
as
\[
N_{p}(f,Y)=\operatorname*{rank}S_{p}(f,Y).
\]
The \textit{Nielsen number for the coincidence problem is denoted by}
$N_{p}(f,g).$
\end{definition}

\begin{corollary}
Suppose $f\sim g.$ Then
\[
N_{\ast}(f)\leq\operatorname*{rank}\,\Omega_{\ast}(g^{-1}(Y)).
\]

\end{corollary}

Clearly $N_{0}(f)$ is equal to the classical Nielsen number and provides a
lower estimate of the number of path components of $f^{-1}(Y).$

It is easy to verify that this theory is still valid if the oriented bordism
$\Omega_{\ast}$ is replaced with the unoriented bordism, or the framed bordism
(see examples in Section \ref{Examples}), or bordism with coefficients. In
fact, a similar theory for an arbitrary homology theory is valid because every
homology theory can be constructed as a bordism theory with respect to
manifolds with singularities \cite{BRS}.

\section{Naturality of $S_{\ast}(f).\label{Natur}$}

Under conditions of Definition \ref{def k*}, the homomorphism $k_{\ast
}:S_{\ast}(f)\rightarrow S_{\ast}(g)$ can be defined as a restriction of
$k_{\ast}^{\prime}$ and the analogues of Propositions \ref{k*} -
\ref{identity} hold. We simplify the situation in comparison to Section
\ref{NClass} by assuming that $Z$ and $Y\subset Z$ are fixed.

\begin{definition}
Suppose we have another preimage problem $g:U\rightarrow Z$ connected to the
first by map $k:X\rightarrow U$ such that $gk=f.$ Then the
\textit{homomorphism induced by }$k,$%
\[
k_{\ast}:S_{\ast}(f)\rightarrow S_{\ast}(g),
\]
is defined as the restriction of $k_{\ast}^{\prime}:S_{\ast}^{\prime
}(f)\rightarrow S_{\ast}^{\prime}(g)$ on $S_{\ast}(f)\subset S_{\ast}^{\prime
}(f).$
\end{definition}

\begin{proposition}
$k_{\ast}$ is well defined.
\end{proposition}

\begin{proof}
For convenience let $f=f_{0},g=g_{0},k=k_{0}.$ Suppose $G$ is a homotopy
between $g_{0}$ and $g_{1},$ $K$ between $k_{0}$ and $k_{1}.$ Let $F=GK,$ then
$F$ is a homotopy between $f_{0}$ and $f_{1}.$ Let $L(t,x)=(t,K(t,x)).$ Then
we have a commutative diagram:%
\[%
\begin{array}
[c]{ccccc}%
U & ^{\underrightarrow{\quad j_{1}\quad}} & I\times U & ^{\underleftarrow
{\quad j_{0}\quad\ }} & U\\
\uparrow^{k_{1}} &  & \uparrow^{L} &  & \uparrow^{k_{0}}\\
X & ^{\underrightarrow{\quad i_{1}\quad}} & I\times X & ^{\underleftarrow
{\quad i_{0}\quad\ }} & X,
\end{array}
\]
where $i_{s}:X\rightarrow\{s\}\times X\rightarrow I\times X$ and
$j_{s}:U\rightarrow\{s\}\times U\rightarrow I\times U,s=0,1,$ are the
inclusions. Further, if we add a vertex $Z$ to this diagram we have a
commutative pyramid with the other edges provided by $f_{0},f_{1},g_{0}%
,g_{1},G,F.$ Then by naturality of the map induced on $S_{\ast}^{\prime}$
(Proposition \ref{composition}) we have another commutative diagram:%
\[%
\begin{array}
[c]{ccccc}%
S_{\ast}^{\prime}(g_{1}) & ^{\underrightarrow{\quad j_{1\ast}^{\prime}\quad}}
& S_{\ast}^{\prime}(G) & ^{\underleftarrow{\quad j_{0\ast}^{\prime}\quad\ }} &
S_{\ast}^{\prime}(g_{0})\\
\uparrow^{k_{1\ast}^{\prime}} &  & \uparrow^{L_{\ast}^{\prime}} &  &
\uparrow^{k_{0\ast}^{\prime}}\\
S_{\ast}^{\prime}(f_{1}) & ^{\underrightarrow{\quad i_{1\ast}^{\prime}\quad}}
& S_{\ast}^{\prime}(F) & ^{\underleftarrow{\quad i_{0\ast}^{\prime}\quad\ }} &
S_{\ast}^{\prime}(f_{0}).
\end{array}
\]
Here the horizontal arrows are injective (Proposition \ref{k-monic}).
Therefore the restriction $k_{0\ast}^{\prime}=k_{1\ast}^{\prime}=L_{\ast
}^{\prime}:M_{\ast}^{F}\rightarrow M_{\ast}^{G}$ is well defined. This
conclusion is true for all $G,K,$ so that the restriction $k_{0\ast}^{\prime
}:\cap_{F=GK}M_{\ast}^{F}\rightarrow\cap_{G}M_{\ast}^{G}$ is well defined.
Since $S_{\ast}(f)$ is a subset of the former and the latter is $S_{\ast}(g)$,
the statement follows.
\end{proof}

\begin{proposition}
Suppose the following diagram for three preimage problems commutes:
\[%
\begin{array}
[c]{ccccc}
&  & Z &  & \\
& \overset{f}{\nearrow} & \uparrow^{f^{\prime}} & \nwarrow^{f^{\prime\prime}}
& \\
X & ^{\underrightarrow{\quad k\quad}} & X^{\prime} & ^{\underrightarrow{\quad
j\quad}} & X^{\prime\prime}%
\end{array}
\]
Then $j_{\ast}k_{\ast}=(jk)_{\ast}:S_{\ast}(f)\rightarrow S_{\ast}%
(f^{\prime\prime}).$
\end{proposition}

\begin{proposition}
$(Id_{X})_{\ast}=Id_{S_{\ast}(f)}.$
\end{proposition}

\begin{corollary}
Given $Z,Y\subset Z.$ If $\mathcal{P}(Z,Y)$ is the category of preimage
problems as pairs $(X,f),$ $f:X\rightarrow Z,$ with morphisms as maps
$k:X\rightarrow U$ satisfying $gk=f,$ then $S_{\ast}$ is a functor from
$\mathcal{P}(Z,Y)$ to $\mathbf{Ab}_{\ast}$ (cf. \cite[Chapter 3]{Jiang}).
\end{corollary}

\begin{corollary}
If $k$ is a homotopy invariance, $gk=f,$ then $S_{\ast}(f)=S_{\ast}(g),$ i.e.,
the preimage theory $f:X\rightarrow Z\supset Y$ is \textquotedblleft homotopy
invariant\textquotedblright\ (cf. \cite{GNS}) with respect to $X.$
\end{corollary}

\section{The Bordism Index as a Homomorphism on $S_{\ast}^{\prime}%
(f)$.\label{BordInd}}

In the classical Nielsen theory, the coincidence index provides an algebraic
count of coincidence points. It satisfies the usual properties: (1) Homotopy
Invariance: the index is invariant under homotopies of $f,g$; (2) Additivity:
The index over a union of disjoint sets is equal to the sum of the indices
over these sets; (3) Existence of Coincidences: if the index is nonzero then
there is a coincidence; (4) Normalization: the index is equal to the Lefschetz
number; (5) Removability: if the index is zero then a coincidence can be
(locally or globally) removed by a homotopy. From the point of view of our
approach the additivity property means that we associate an integer to every
$0$-class, i.e., we have a homomorphism $S_{0}^{\prime}(f)\rightarrow
\mathbf{Z=}\Omega_{0}(Y)$.

\begin{definition}
We define \textit{the index} of $[s]_{N}\in S_{q}^{\prime}(f)$, where
$s\in\Omega_{q}(C),$ as
\[
I_{f}(s)=f_{\ast}(s),
\]
where $f_{\ast}:\Omega_{q}(C)\rightarrow\Omega_{q}(Y).$
\end{definition}

\begin{proposition}
The index is well defined as a homomorphism
\[
I_{f}:S_{\ast}^{\prime}(f)\rightarrow\Omega_{\ast}(Y).
\]

\end{proposition}

\begin{proof}
Suppose $s\sim_{N}\emptyset,$ then by definition, $s\sim_{b}\emptyset$ in $X$
via some $H$ and $fs\sim_{b}\emptyset$ in $Z$ via some $G$ homotopic to $fH.$
Therefore $f_{\ast}(s)=0$ in $\Omega_{q}(Y),$ so $I_{f}([s]_{N})=0.$
\end{proof}

Of course, $I_{f}(z)\neq0\Longrightarrow z\neq0.$

Suppose $F:I\times X\rightarrow Z$ is a homotopy. As before, let $f_{t}%
(\cdot)=F(t,\cdot):X\rightarrow Z,$ and let $i_{t}:X\rightarrow\{t\}\times
X\rightarrow I\times X$ be the inclusions.

\begin{theorem}
Suppose $z_{0}\in S_{q}^{\prime}(f_{0})$ and $i_{0\ast}^{\prime}(z_{0})=z\in
S_{q}^{\prime}(F).$ Then $I_{f_{0}}(z_{0})=I_{F}(z).$
\end{theorem}

\begin{proof}
Let $C=F^{-1}(Y),C_{0}=f_{0}^{-1}(Y).$ Now if $z=[s]_{N},s\in\Omega_{q}(C)$
and $z_{0}=[s_{0}]_{N},s_{0}\in\Omega_{q}(C_{0}),$ then $i_{0\ast}(s_{0}%
)\sim_{N}s.$ Therefore, $I_{F}(z)=F_{\ast}(s)=F_{\ast}i_{0\ast}(s_{0}%
)=f_{0\ast}(s_{0})=I_{f_{0}}(z_{0}).$
\end{proof}

\begin{corollary}
If $z_{0}\in S_{q}^{\prime}(f_{0}),$ $z_{1}\in S_{q}^{\prime}(f_{1})$ are
$F$-related then $I_{f_{0}}(z_{0})=I_{f_{1}}(z_{1}).$ Thus the index $I_{f}$
is preserved under homotopies.
\end{corollary}

In the classical theory Nielsen classes are sets and the algebraically
essential classes are the ones with nonzero index. Similarly, we call $z\in
S_{q}^{\prime}(f)$ \textit{algebraically essential} of $I_{f}(z)\neq0.$

\begin{corollary}
Every algebraically essential class is topologically essential, i.e., $z$
cannot be reduced by a homotopy to the zero $p$-class, and, therefore, $z$
cannot be \textquotedblleft removed\textquotedblright\ by a homotopy.
\end{corollary}

We define the \textit{group of algebraically essential Nielsen classes }as
\[
S_{\ast}^{a}(f,Y)=S_{\ast}^{\prime}(f,Y)/\ker I_{f}.
\]

Suppose we have another preimage problem $g:U\rightarrow Z$ connected to the
first by map $k:X\rightarrow U$ such that $gk=f.$ Then just like in the
previous section we define the \textit{homomorphism induced by }$k,$%
\[
k_{\ast}^{a}:S_{\ast}^{a}(f)\rightarrow S_{\ast}^{a}(g),
\]
as a restriction of $k_{\ast}^{\prime}.$ Moreover similar properties are
satisfied. Thus we have

\begin{corollary}
Given $Z,Y\subset Z.$ If $\mathcal{P}(Z,Y)$ is the category of preimage
problems as pairs $(X,f),$ $f:X\rightarrow Z,$ with morphisms as maps
$k:X\rightarrow U$ satisfying $gk=f,$ then $S_{\ast}^{a}$ is a functor from
$\mathcal{P}(Z,Y)$ to $\mathbf{Ab}_{\ast}$.
\end{corollary}

\section{The Index of an Isolated Set of Preimages.\label{Index}}

From now on we assume that $X,Y,Z$ are smooth manifolds, $Y$ is a submanifold
of $Z.$

Suppose $P\subset C$ is an isolated set of preimages. Let $U\subset V$ be
neighborhoods of $P$ in $X$ such that $\overline{U}\subset Int\,V$ and $V\cap
C=P.$ In the classical Nielsen theory, the index $Ind(f,P)$ of $P$ is defined
as the image of the generator $z$ of $H_{n+m}(X)\simeq\mathbf{Z}$ under the
composition
\[
H_{n+m}(X)^{\underrightarrow{\qquad}}H_{n+m}(X,X\backslash U)^{\underleftarrow
{\quad\simeq\quad\ }}H_{n+m}(V,V\backslash U)^{\underrightarrow{\quad f_{\ast
}\quad}}H_{n+m}(Z,Z\backslash Y).
\]
Under the restriction $\dim X+\dim Y=\dim Z,$ we are in the classical
situation: each class $a$ is an isolated set of preimages $A$ and the index
$a$ is defined as the index of $A$.

In case of a nonzero codimension we can have $H_{n+m}(X)\neq\mathbf{Z,}$
therefore it makes sense to replace in the above definition the generator $z$
with an arbitrary element of $H_{\ast}(X)$. This turns the index into a graded
homomorphism $H_{\ast}(X)\rightarrow H_{\ast}(Z,Z\backslash Y)$ (which is
equal to the Lefschetz homomorphism \cite{Sav1} for the Coincidence Problem).
This generality is justified by a number of examples in \cite{GNO},
\cite{Sav1} that show that in order to detect coincidences in a nonzero
codimension one may need to take into account all parts of this homomorphism.

A fixed point index with respect to generalized cohomology was considered by
Dold \cite{Dold1}. Another example is \cite{Mukh} where the coincidence index
is computed in term of cobordism. In addition, we will see in Section
\ref{Jez} that for a nonzero codimension the index expressed in terms of
singular homology may be inadequate for removability (some algebraically
inessential classes are essential). Therefore under the above restrictions the
singular homology $H_{\ast}$ should be replaced with a generalized homology
$h_{\ast}$.

\begin{definition}
Suppose $W$ is a neighborhood of $Y$ in $Z.$ The \textit{index} $Ind_{f}%
(P;h_{\ast}),$ or simply $Ind_{f}(P),$ of the set $P$ with respect to
$h_{\ast}$ is the following homomorphism
\[
h_{\ast}(X,X\backslash U)^{\underleftarrow{\quad\simeq\quad\ }}h_{\ast
}(V,V\backslash U)^{\underrightarrow{\quad f_{\ast}\quad}}h_{\ast
}(W,W\backslash Y).
\]

\end{definition}

The index does not depend on the choice of $U$, see \cite[p. 189]{Vick}. The
next theorem is proven similarly to Lemmas 7.1, 7.2, 7.4 in \cite[p.
190-191]{Vick} respectively.

\begin{theorem}
I. (Additivity) If $P=P_{1}\cup...\cup P_{s}$ is a disjoint union of subsets
of $P$ such that $P_{i}=P\cap f^{-1}(Y)$ is compact and $f_{i}=f|_{P_{i}}$ for
each $i=1,...,s,$ then $Ind_{f}(P)=\sum_{i=1}^{s}Ind_{f_{i}}(P_{i}).$

II. (Existence of Preimages) If $Ind_{f}(P)\neq0$ then $P\neq\emptyset.$

III. (Homotopy Invariance) Suppose $f_{t}:X\rightarrow Z,0\leq t\leq1,$ is a
homotopy and $D=\cup_{t}f_{t}^{-1}(Y)$ is a compact subset of $X.$ Then
\[
Ind_{f_{0}}(P)=Ind_{f_{1}}(P).
\]

\end{theorem}

Now we show how the indices $I_{f}$ and $Ind_{f}$ are related to each other$.$

\begin{definition}
Let $R_{p}(f)$ be the set of all elements of $M_{p}(C)$ represented as sums of
finite collections of connected singular manifolds such that the sum of any
subcollection is not Nielsen equivalent to the empty set. Then \textit{the
image of} $z\in S_{p}^{\prime}(f)$ is defined as
\[
\operatorname{Im}z=\cup\{\operatorname{Im}s:s\in R_{p}(f),s\in z\}.
\]

\end{definition}

In particular, if $P\in\Omega_{0}(M)$ is a subset of $C$ and $z=[P]_{N},$ then
$\operatorname{Im}z=P.$

Suppose $f$ is transversal to $Y$. Then $C$ is a $r$-submanifold of $X$.

\begin{proposition}
$\operatorname{Im}z$ is an isolated subset of $C$ and therefore an
$r$-submanifold of $X.$
\end{proposition}

\begin{proposition}
If $z=0\in S_{\ast}^{\prime}(f)$ then $\operatorname{Im}z=\emptyset.$
\end{proposition}

Clearly if $Ind_{f}(\operatorname{Im}z)\neq0$ then $\operatorname{Im}%
z\neq\emptyset.$ However this does not imply that $z$ is essential. The case
of $p=0$ is an exception. For convenience we restate the following familiar result.

\begin{proposition}
\label{alg-ess-0}If $P\in S_{0}^{\prime}(f)$ and $Ind_{f}(P;h_{\ast})\neq0,$
where $h_{\ast}$ is an arbitrary homology theory, then $P$ is essential.
\end{proposition}

The relation between the essentiality of the class and its index is more
subtle when $p>0.$

If $T$ is a tubular neighborhood of a submanifold $M,$ \ then $\varphi
_{M}:\Omega_{q+k}(T,T\backslash M)\rightarrow\Omega_{q}(M)$ is the Thom
isomorphism \cite[p. 309]{Switzer}, \cite[p. 321]{Dold}.

Suppose $z\in S_{p}^{\prime}(f),z=[s],$ where $s\in\Omega_{p}(M).$ Let
$P=\operatorname{Im}z,$ then it is an $r$-submanifold of $C.$ Let $T$ and
$T^{\prime}$ be tubular neighborhoods of $C$ and $P$ respectively such that
$T^{\prime}$ is an isolated subset of $T.$ Then the inclusion $i:T^{\prime
}\rightarrow T$ is a bundle map. Suppose $s=i_{\ast}(s^{\prime})$ for some
$s^{\prime}\in\Omega_{p}(P).$ From the naturality of the Thom isomorphism we
have the commutativity of the following diagram:%
\[%
\begin{tabular}
[c]{cccccccc}%
$\Omega_{p}(P)$ & $^{\underrightarrow{\ \ i_{\ast}\ \ }}$ & $\Omega_{p}(C)$ &
$^{\underrightarrow{\ \ f_{\ast}\ \ }}$ & $\Omega_{p}(Y)$ &  &  & \\
$\uparrow^{\varphi_{P}}$ &  & $\uparrow^{\varphi_{C}}$ &  & $\uparrow
^{\varphi_{Y}}$ &  &  & \\
$\Omega_{p+k}(T^{\prime},T^{\prime}\backslash P)$ & $^{\underrightarrow
{\ \ i_{\ast}\ \ }}$ & $\Omega_{p+k}(T,T\backslash C)$ & $^{\underrightarrow
{\ \ f_{\ast}\ \ }}$ & $\Omega_{p+k}(W,W\backslash Y),$ &  &  &
\end{tabular}
\ \
\]
where $W$ is a tubular neighborhood of $Y.$ Then $I_{f}(z)=f_{\ast}%
(s)=f_{\ast}i_{\ast}(s^{\prime}).$ Therefore
\[
\varphi_{Y}^{-1}I_{f}(z)=\varphi_{Y}^{-1}f_{\ast}i_{\ast}(s^{\prime})=f_{\ast
}i_{\ast}\varphi_{P}^{-1}(s^{\prime})=Ind_{f}(P;\Omega_{\ast})(\varphi
_{P}^{-1}(s^{\prime})).
\]
Thus we have proven

\begin{theorem}
\label{bi}For each $z\in S_{p}^{\prime}(f)$, $I_{f}(z)=\varphi_{Y}%
Ind_{f}(P;\Omega_{\ast})\varphi_{P}^{-1}(s^{\prime}),$ where $z=[s]$ and
$s=i_{\ast}(s^{\prime}).$
\end{theorem}

The right hand side can be used for an alternative definition of an
algebraically essential class.

\begin{corollary}
If $P\in\Omega_{r}(M),z=[P]$, then $Ind_{f}(P;\Omega_{\ast})(O_{P}%
)=\varphi_{Y}^{-1}I_{f}(z),$ where $O_{P}\in\Omega_{n+m}(T^{\prime},T^{\prime
}\backslash P)$ is the fundamental class of $P.$
\end{corollary}

Moreover if $Ind_{f}(P;H_{\ast})\neq0$ then $P\in S_{0}^{\prime}(f)$ is
essential. Thus for $r=0$ we recover the traditional definition of an
algebraically essential class.

\section{Some Examples.\label{Examples}}

Nielsen numbers are hard to compute. Nielsen groups and higher order Nielsen
numbers are no different. Below we consider some special cases when the
computation is feasible.

Just as before suppose $X=Z=\mathbf{S}^{2},$ $Y$ is the equator of $Z,$ $f$ a
map of degree $2$ such that $C=f^{-1}(Y)$ is the union of two circles $C_{1}$
and $C_{2}$ around the poles. But there is only one generator of
$S_{1}^{\prime}(f),$ $C=C_{1}\cup C_{2}\mathbf{.}$ Also $Ind_{f}(C)\neq0$.
Hence $N_{0}(f)=1.$ This is in fact a \textquotedblleft
sharp\textquotedblright\ estimate of the number of components of $C$ (Wecken
Property for codimension $r=1$) because $f$ is homotopic to the suspension,
$g,$ of the degree $2$ map of the equator, so that $g^{-1}(Y)$ is a circle.
The same conclusion applies to $X=Z=\mathbf{S}^{n},Y=\mathbf{S}^{n-1},$
codimension $r=n-1,n\geq2.$

For more examples of this nature see \cite{RW} and Theorem 1.2 and Section 12
in \cite{DG}.

In \cite{Sav2} we showed that the cohomology coincidence index $I_{fg}^{A}$ is
the only obstruction to removability of an isolated subset $A$ of the
coincidence set if any of the three following conditions is satisfied: (1) $M$
is a surface; (2) the fiber of $g$ is acyclic; or (2) the fiber of $g$ is an
$m$-sphere for $m=4,5,12$ and $n$ large. Of course if the homology index
$Ind_{(f,g)}(A;H_{\ast})$ is trivial then so is $I_{fg}^{A}.$ Therefore under
these restrictions an algebraically inessential class can be removed.

\begin{proposition}
Suppose that if two singular $q$-manifolds in $C$ are bordant in $X$ then they
are Nielsen equivalent. Then $S_{q}^{\prime}(f)\simeq j_{\ast}\Omega_{q}(C),$
where $j:C\rightarrow X$ is the inclusion.
\end{proposition}

\begin{corollary}
If under the conditions of the proposition $\Omega_{q}(X)=0,$ then
$S_{q}(f)=0.$
\end{corollary}

The condition of this proposition is satisfied if we simply assume that $f$ is
homotopic to $f^{\prime}$ with $f^{\prime}(X)\subset Y.$ For the Coincidence
Problem this result takes the following form.

\begin{theorem}
If $f,g:X\rightarrow Z$ are homotopic then $S_{q}^{\prime}(f,g)\simeq j_{\ast
}\Omega_{q}(Coin(f,g)).$ Moreover, $N_{q}(f,g)\leq\operatorname*{rank}%
\Omega_{q}(X).$
\end{theorem}

Suppose a parametrized dynamical system $F\rightarrow N^{\underrightarrow
{~\ \ \ f,g\ \ \ ~~}}M,$ where $g$ is the bundle projection, is generated by a
control system. Then it is easy to see that $f$ is homotopic to $g$ and the
above theorem can be applied to estimate the \textquotedblleft
size\textquotedblright\ of the set of stationary points of the control system.

\begin{theorem}
Suppose $Y$ is $(q-1)$-connected, $f^{\ast}:H^{q}(Y;\pi_{q}(Y))\rightarrow
H^{q}(X;\pi_{q}(Y))$ is trivial, and $Z$ is $(q+1)$-connected. Then
$S_{q}^{\prime}(f)\simeq j_{\ast}\Omega_{q}(C).$ Moreover, $N_{q}%
(f)\leq\operatorname*{rank}\Omega_{q}(X).$
\end{theorem}

\begin{proof}
Suppose $s_{0},s_{1}\in j_{\ast}\Omega_{q}(C)$ are bordant in $X$ via
$H:W\rightarrow X,$ i.e., $s_{i}:S_{i}\rightarrow X,\partial W=S_{0}\sqcup
S_{1},H|_{S_{i}}=s_{i}.$ Since $Y$ is $(q-1)$-connected and
\[
\delta^{\ast}f^{\ast}(s_{0}\sqcup s_{1})^{\ast}:H^{q}(Y;\pi_{q}(Y))\rightarrow
H^{q+1}(W,S_{0}\sqcup S_{1};\pi_{q}(Y))
\]
is trivial, the classical obstruction theory \cite[p. 497]{Bredon} is applied
to prove that the map $f(s_{0}\sqcup s_{1}):S_{0}\sqcup S_{1}\rightarrow Y$
can be extended to $G:W\rightarrow Y$ . Further, since $Z$ is $(q+1)$%
-connected$,$ $[W,Z]_{relS_{0}\sqcup S_{1}}=0.$ Therefore $G$ and $fH$ are
homotopic relative to $S_{0}\sqcup S_{1}.$ Thus, if two singular $q$-manifolds
in $C$ are bordant in $X,$ then they are Nielsen equivalent. Now the theorem
follows from the above proposition.
\end{proof}

The relation between the homotopy class of a map and the preimage of a point
is direct in the setting of the Pontriagin-Thom construction \cite[p.
196]{Davis}. For the rest of the section we assume that the Nielsen groups
$S_{q}^{\prime}(f),$ $S_{q}(f)$ are computed with respect to the
\textit{framed} bordism, i.e., $S_{q}^{\prime}(f)$ is a quotient group of
$\Omega_{q}^{fr}(C).$

Let $Y=\{p\},p\in Z=\mathbf{S}^{k},$ and $r\leq k-2.$ Then the conditions of
the theorem above are satisfied. Therefore $S_{r}^{\prime}(f)\simeq j_{\ast
}\Omega_{r}^{fr}(C).$ Now, $f$ is homotopic to a map $g$ if and only if
$C=f^{-1}(p)$ is framed bordant to $K_{g}=g^{-1}(p)$ in $X.$ Let $j^{g}%
:K_{g}\rightarrow X$ be the inclusion. Then
\[
S_{r}(f)\simeq%
{\displaystyle\bigcap\limits_{g\sim f}}
j_{\ast}^{g}\Omega_{r}^{fr}(K_{g})=%
{\displaystyle\bigcap\limits_{K_{g}\sim_{b}C}}
j_{\ast}^{g}\Omega_{r}^{fr}(K_{g}).
\]
Thus we have proven the following

\begin{theorem}
\label{PCRT}Suppose $Z=\mathbf{S}^{k}$ and $r=m-k\leq k-2.$ Then
\[
S_{r}(f,\{p\})=%
{\displaystyle\bigcap\limits_{K\sim_{b}C}}
j_{\ast}^{K}\Omega_{r}^{fr}(K),
\]
where $j^{K}:K\rightarrow X$ is the inclusion.
\end{theorem}

In particular for codimension $1,$ $N_{1}(f,\{p\})$ is equal to the number of
circles in $f^{-1}(p)$ not framed bordant to the empty set.

\begin{corollary}
Suppose $Z=\mathbf{S}^{k}$ and $r=m-k\leq k-2.$ Then $S_{r}(f,\{p\})=0$ iff
$f^{-1}(p)\sim_{b}\varnothing.$
\end{corollary}

\begin{proof}
The right hand side in the above theorem contains $C$.
\end{proof}

Nielsen groups can be easily computed for the generators of $[\mathbf{S}%
^{k},\mathbf{S}^{m}],$ see \cite[p. 208]{Davis}.

\section{Wecken Property of Order 1, Codimension $1$.\label{Wecken}}

We say that the preimage problem $f:X\rightarrow Z\supset Y$ satisfies the
\textit{Wecken Property of order} $p$ if $S_{p}(f)$ is \textquotedblleft
realizable\textquotedblright$,$ i.e., there is some $h$ homotopic to $f$ such
that $S_{p}(f)\simeq\Omega_{p}(h^{-1}(Y)).$

Recall that a preimage problem is reduced to the coincidence problem
$f,g:X\rightarrow Y$ by putting $Z=Y\times Y,$ $Y$ the diagonal of $Z,$
$f=(F,G),C=(f,g)^{-1}(Y)=Coin(f,g).$ Also $\dim Z=2\dim Y=2n,$ so $k=n$

Then the definition of Wecken Property is in the obvious extension of the one
above: the pair $(f,g)$ satisfies the \textit{Wecken Property of order} $p$ if
$S_{p}(f,g)\simeq\Omega_{p}(Coin(f^{\prime},g^{\prime}))$ for some $f^{\prime
},g^{\prime}$ homotopic to $f,g.$

Assume that $f,g$ are transversal. Then $C=Coin(f,g)$ is a $1$-submanifold of
$X.$ Suppose $A$ is $1$-submanifold of $C$. Recall condition $\text{(W1)
}A=\partial S,\text{ }$\textit{where} $S$ \textit{is an orientable connected
surface}$,\text{ }f|_{S}\sim g|_{S}\text{ rel }A.$

\begin{proposition}
\label{W1}(W1$^{\prime}$) $A=\partial S,\text{ }$\textit{where} $S$ \textit{is
an orientable (not necessarily connected) surface}$,\text{ }f|_{S}\sim
g|_{S}\text{ rel }A\Longleftrightarrow$ $A\sim_{N}\varnothing,$ i.e., $A$
belongs to the zero $1$-class.
\end{proposition}

\begin{proof}
Since $A$ is the boundary of $S,$ we have $A\sim_{b}\varnothing$ via $S.$
Secondly, $f|_{S}\sim g|_{S}$ rel $A,$ hence $(f,g)|_{S}\sim h$ rel $A$ such
that $\operatorname{Im}h$ lies in the diagonal of $Z=Y\times Y.$ Thus
$A\sim_{N}\varnothing.$
\end{proof}

Consider the following result due to Jezierski \cite[Theorem 3.1]{Jez}. Its
proof is based on his $1$-parameter Whitney Lemma.

\begin{proposition}
[Jezierski]\label{Circles}Suppose $n\geq4$ and $f$ is smooth. Then $g$ is
homotopic to $g^{\prime}$ such that the pair $(f,g^{\prime})$ is transversal
and each Nielsen class (or even an isolated subset of a Nielsen class
\cite{Jez-p}) of $(f,g)$ is a circle.
\end{proposition}

We use this result to prove the following Wecken-type result for codimension
$1.$

\begin{theorem}
\label{WP0}Suppose $n\geq4$ and $f$ is smooth. Then the coincidence problem
$f,g:X\rightarrow Y$ satisfies the \textit{Wecken Property of order} $1$;
specifically, $g$ is homotopic to $g^{\prime}$ such that%
\[
S_{1}(f,g)\simeq\Omega_{1}(Coin(f,g^{\prime})).
\]
Moreover, $N_{1}(f,g)=\operatorname*{rank}\Omega_{1}(Coin(f,g^{\prime}))$ is
equal to the number of circles in $Coin(f,g^{\prime})$ not satisfying (W1).
\end{theorem}

\begin{proof}
Suppose, according to the above proposition, that all $0$-classes are circles,
$A_{1},...,A_{s}.$ Suppose also that $A_{1},...,A_{t}$ satisfy condition (W1)
and the rest do not. Let's view $A_{1},...,A_{t}$ as singular $1$-manifolds.
Then, first, $A_{i}\sim_{N}\varnothing$ for $i=1,...,t$ according to
Proposition \ref{W1}. Hence for these $i,$ $A_{i}\in0\in S_{1}^{\prime
}(f,g^{\prime}),$ so they don't concern us. Now, suppose $A_{i}\sim_{N}A_{j}$
for some $i>j>t$ via some surface $H.$ If $A_{i}$ and $A_{j}$ were subsets of
different components of $H$ then each would satisfy condition (W1). Therefore
$H$ can be assumed connected. But then $A_{i}\cup A_{j}$ satisfies condition
(W1) and, moreover, every pair of points $x\in A_{i},y\in A_{j}$ is Nielsen
equivalent. Therefore by Proposition \ref{Circles} $A_{i}\cup A_{j}$ can be
further reduced to a single circle. Hence we can assume that each
$A_{i},i=t+1,...,s,$ belongs to a different nonzero $1$-class. Thus the
generators of $S_{1}^{\prime}(f,g^{\prime})$ are $[A_{i}]_{N},i=t+1,...,s.$
Now the fact that each of these classes is essential follows from the homotopy
invariance of (W1) \cite[Lemma 4.2]{Jez}.
\end{proof}

A similar result for the Root Problem is easy to prove.

\begin{theorem}
Suppose $Z=\mathbf{S}^{k}$ and $1=m-k\leq k-2.$ Then the root problem
$f:X\rightarrow Z\ni p$ satisfies the Wecken Property of order $1$ (with
respect to framed bordism).
\end{theorem}

\begin{proof}
Just as above assume that $C=f^{-1}(p)$ is the disjoint union of circles such
that $A_{1},...,A_{t}$ are framed bordant to the empty set and $A_{t+1}%
,...,A_{s}$ are not. Then $C$ is framed bordant to $K=A_{t+1}\cup...\cup
A_{s}.$ Finally $S_{1}(f,\{p\})=\Omega_{1}^{fr}(K)$ by Theorem \ref{PCRT}.
\end{proof}

\section{Wecken Property of Order $0,$ Codimension $1.\label{Jez}$}

Suppose $A$ is a $1$-submanifold of $C$. Recall condition (W2)\textit{ the PT
map is trivial. }The proposition below explains why the PT map should be
understood as the coincidence index.

\begin{proposition}
(W2) $\Longleftrightarrow Ind_{(f,g)}(A;\pi_{\ast}^{S})=0$ (i.e., $A$ is
algebraically inessential with respect to $\pi_{\ast}^{S}).$
\end{proposition}

\begin{proof}
Let $U\subset T$ be tubular neighborhoods $A.$ We state (W2) as follows:
\[
PT:\mathbf{S}^{n+1}\rightarrow T/\partial T^{\underrightarrow
{~\ \ \ f-g~\ \ ~~}}\mathbf{R}^{n}/(\mathbf{R}^{n}\backslash0)\simeq
\mathbf{S}^{n}\text{ is trivial.}%
\]
Since $n\geq4,$ this is equivalent to the following:
\[
PT_{\ast}:\pi_{\ast}^{S}(\mathbf{S}^{n+1})^{\underrightarrow{~\ \ \ ~~}}%
\pi_{\ast}^{S}(\mathbf{S}^{n})\text{ is trivial.}%
\]
Consider the commutative diagram:%
\[%
\begin{array}
[c]{ccccccc}
&  &  & \nearrow & T/\partial T & ^{\underrightarrow{~\ \ \ f-g~\ \ ~~}} &
(\mathbf{R}^{n},\mathbf{R}^{n}\backslash0)\\
\mathbf{S}^{n+1} & \rightarrow & \mathbf{S}^{n+1}/(\mathbf{S}^{n+1}\backslash
U) &  & \downarrow^{\simeq} &  & ^{\simeq}\uparrow^{d}\\
&  &  & \searrow & (T,T\backslash U) & ^{\underrightarrow{~\ \ \ (f,g)~\ \ ~~}%
} & (\mathbf{R}^{n}\times\mathbf{R}^{n},\mathbf{R}^{n}\times\mathbf{R}%
^{n}\backslash\Delta),
\end{array}
\]
where $\Delta$ is the diagonal and $d(x,y)=x-y.$ Now if we apply the stable
homotopy functor $\pi_{\ast}^{S}$ to the diagram, we have $PT_{\ast}$ in the
upper path and the index of $A$ with respect to $h_{\ast}=\pi_{\ast}^{S}$ in
the lower. But $d$ is a homotopy equivalence \cite[Lemma VII.4.13, p.
200]{Dold}, and the statement follows.
\end{proof}

Observe that the stable homotopy index $Ind_{(f,g)}(A;\pi_{\ast}^{S})$ is
better at detecting essential classes than the traditional index with respect
to singular homology. In fact, the latter would not work in the above argument
as $\pi_{n+1}^{S}(\mathbf{S}^{n})=\mathbf{Z}_{2}$ cannot be replaced with
$H_{n+1}(\mathbf{S}^{n})=0$. Secondly, all the Nielsen numbers of higher order
in Section \ref{Examples} would be zero if computed with respect to singular homology.

Recall Jezierski's Wecken type theorem \cite[Theorem 5.3]{Jez}.

\begin{proposition}
[Jezierski]\label{JezThm} Let $f,g:X\rightarrow Y$ be an admissible map
between open subsets of $\mathbf{R}^{n+1},$ $\mathbf{R}^{n}$ respectively,
$n\geq4.$ Then there are maps $f^{\prime},g^{\prime}$ compactly homotopic to
$f,g$ respectively such that the Nielsen classes satisfying (W1) and (W2)
disappear and the remaining ones become circles.
\end{proposition}

Suppose we are left with the circles $A_{1},A_{2},...,A_{t}$ each satisfying
condition (W1) but not (W2), and $A_{t+1},A_{t+2},...,A_{s}$ satisfying (W2)
but not (W1). Then each $A_{1},A_{2},...,A_{t}$ is an (algebraically)
essential $0$-class (Theorem \ref{alg-ess-0}). Also each $A_{t+1}%
,A_{t+2},...,A_{s}$ is an essential $1$-class (Theorem \ref{WP0}), therefore
an essential $0$-class as well.\ Thus we have proven the following

\begin{theorem}
\label{JezCor}Suppose $X,Y$ are open subsets of $\mathbf{R}^{n+1},$
$\mathbf{R}^{n}$ respectively, $n\geq4.$ Then there are maps $f^{\prime
},g^{\prime}$ compactly homotopic to $f,g$ respectively such that
$Coin(f^{\prime},g^{\prime})$ has exactly $N_{0}(f,g)$ path components, i.e.,
the coincidence problem $f,g:X\rightarrow Y$ satisfies the Wecken Property of
order $0.$
\end{theorem}

A result of this type is proven by Gon\c{c}alves and Wong \cite[Theorem 4
(iii)]{GW} for the Root Theory and an arbitrary codimension. In the
terminology of the present paper their theorem reads as follows: If $X,Z$ are
nilmanifolds, $p\in Z,$ then there is $g$ homotopic to $f$ such that
$g^{-1}(p)$ has exactly $N_{0}(f,\{p\})$ components.

Another codimension $1$ Wecken type theorem is given by Dimovski \cite{Dim}
for the parametrized fixed point problem: $F:I\times Y\rightarrow Y.$ He
defines two independent indices of a Nielsen class $V,$ $ind_{1}(F,V)$ and
$ind_{2}(F,V),$ corresponding to conditions (W1), (W2), and then defines a
Nielsen number $N(F)$ as the number of Nielsen classes with either
$ind_{1}(F,V)\neq0$ or $ind_{2}(F,V)\neq0.$ His Theorem 4.4 (4) reads: If $F$
is homotopic to $H$ such that $H$ has only isolated circles of fixed points
and isolated fixed points then the number of fixed points classes of $H$ is
bigger than or equal to $N(F).$ However this does not mean that $N(F)$ is a
lower bound of the number of components of the fixed point set because an
examination of the proof of this theorem (Theorems 4.1, 4.2) reveals that only
\textit{local} homotopies, i.e., ones constant outside a neighborhood of the
given class, are allowed. This is the reason why there is a direct
correspondence between Nielsen classes of two homotopic maps and there is no
need for such a construction as the one in Section \ref{TopEss} of the present
paper. Also $N(F)$ can be larger than the estimate provided in the above
theorem - Jezierski \cite[Example 6.4]{Jez} gives an example of a Nielsen
class that can be removed by a global homotopy but not by a local one.

\end{document}